\documentstyle[12pt]{article}
\newcommand{\R}{\mbox{I$\!$R}}

\newcommand{\qed}{{\hfill {$\rlap{$\sqcap$}\sqcup$}}\\[0.2in]\hspace*{0.5in}}
\newcommand{\qedwh}{{\hfill {$\rlap{$\sqcap$}\sqcup$}}\\[0.2in]}
\newcommand{\bk}{\\[0.03in] \hspace*{0.5in} }
\newcommand{\btd}{\bigtriangledown}

\newcommand{\mfor}{\ \ \ \ {\mbox{for}} \ \ }

\begin{document}

\vspace*{0.02in}
\begin{center} 
{\LARGE {\bf  Growth Estimates on Positive Solutions}}
\medskip \\ {\LARGE {\bf  of the Equation $\Delta u + K u^{{n + 2}\over {n - 2}} = 0$ in ${\R}^n$}}
\medskip \medskip \medskip  \medskip\\ 
{\large {Man Chun {\Large L}EUNG\footnote{
Department of
Mathematics,  National University of Singapore, 2 Science Drive 2,\\
Singapore 117543, Republic of Singapore; \ \
{\tt matlmc@math.nus.edu.sg} \\
$\bullet$ KEY WORDS: positive solution, conformal scalar curvature equation, growth estimate.\\
$\bullet$ 2000 AMS MS CLASSIFICATIONS: Primary 35J60\,; Secondary 53C21.\\
$\bullet$ RUNNING TITLE: Growth Estimates on Positive Solutions...
 }}}
\end{center}

\vspace*{0.15in} 

\begin{abstract}We construct unbounded positive
$C^2$-solutions of the equation $\Delta u + K u^{(n + 2)/(n - 2)} = 0$ in 
${\R}^n$ (equipped with Euclidean metric $g_o$)  such that $K$ is bounded between 
two positive numbers in ${\R}^n$, the conformal metric $g = u^{4/(n - 2)} \,g_o$ is complete,
and the volume growth of $g$ can be arbitrarily fast or reasonably slow according to the
constructions. By imposing natural conditions on $u$, we obtain growth estimate on the $L^{2n/(n - 2)}$-norm of the solution and show that it has slow decay.
\end{abstract}

\vspace*{0.5in}

{\bf \Large {\bf 1. \ \ Introduction}}

\vspace*{0.3in}

In this article we derive $L^p$-estimates on positive solutions of the conformal scalar
curvature equation
$$
\Delta u +  K u^{{n + 2}\over {n - 2}} = 0 \ \ \ \ {\mbox{in}} \ \ {\R}^n, \leqno (1.1)
$$
where $n \ge 3$ is an integer, $\Delta$ the standard Laplacian on ${\R}^n,$ 
 $K$ a smooth function. Equation (1.1) relates the scalar curvature
of  the conformal metric $g = u^{4/(n - 2)} \,g_o$ to $4K(n - 1)/(n - 2)$, 
where $g_o$ is  Euclidean metric \cite{Leung.1}. It is assumed
throughout this note that 
$$
0 < a^2 \le K (x) \le b^2 \mfor {\mbox{large}} \ \ |x|  \leqno (1.2)
$$
and for some positive constants $a$ and $b$. The following estimates are
known for any positive smooth solution $u$ of equation (1.1) with condition (1.2).
$$
\int_{S^{n - 1}} u (r, \theta) \,
d
\theta \le  C_1\,  r^{{2 - n}\over 2}\,, \leqno (1.3)
$$
$$
\int_{B_o (r)} u^{{n + 2}\over
{n - 2}} (x) \, dx\le  C_2 \, r^{{n - 2}\over 2} \leqno (1.4)
$$
for large $r$ and for some positive constants $C_1$ and $C_2$ depending on $u$ (see, for
example,
\cite{Leung.5}). Here $B_o (r)$ is the ball
with center at the origin and radius $r\,,$ and $S^{n - 1}$ is the unit sphere in 
${\R}^n$. We seek to obtain higher order estimates of the forms 
$$
\int_{S^{n - 1}} u^p (r, \theta) \, d \theta \le C_3\,  r^{ (2 - n) p/2},  \ \ \ \ p >
1\,;\leqno (1.5)
$$
\[ (1.6) \ \ \ \ \ \ \ \ \ \int_{B_o (r)} u^q (x) \, dx \le  \left\{ \begin{array}   
{r@{\quad \mbox{if}\quad}l}
C_4 \,  r^{ n - (n - 2) q/2} \  & \ \  q > {{n + 2}\over {n - 2}} \,, \ q\not= {{2n}\over
{n - 2}}\,; \
\
\
\
\
\
\
\
\
\  \ \ \ \ \ \ \ \ \ \\ 
 C_5 \, \ln r    \ \ \ \ & \ \  q = {{2n}\over {n - 2}}\,, 
\end{array} \right.  \]
for large $r$, where  $C_3\,, \ C_4$ and $C_5$ are positive constants. 
The above estimates  are based on the {\it slow decay} of $u$, that is,
$$
u (x) \le C_6 \, |x|^{(2 - n)/2}  \mfor {\mbox{large}} \ \ |x|\,, \leqno (1.7)
$$
where $C_6$ is a positive constant.\bk 
Our first
observation is that, in general, (1.5), (1.6) or (1.7) do not hold.  Taliaferro
\cite{Taliaferro.1} shows that positive solution of (1.1) outside a ball in ${\R}^n$ 
with condition (1.2) may not have slow decay. We modify the construction in
\cite{Taliaferro.1} to obtain positive
$C^2$-solutions of (1.1) in 
${\R}^n$  with $K$ bounded between two positive numbers in ${\R}^n$, 
such that the conformal metric $g = u^{4/(n - 2)} \,g_o$ is complete and the volume growth of $({\R}^n, g)$ can
be arbitrarily fast or reasonably slow with respect to the constructions. This suggests that the geometric structure of complete manifolds of bounded positive scalar curvature could be very complicated (cf. \cite{Gromov.1}).\bk
It is observed in
\cite{Cheung-Leung.1} that if estimate (1.5) holds for some number $p \ge 2n/(n - 2)$, then $u$ has
slow decay and hence (1.5) and (1.6) hold for all $p,\, q > 1$. 
The integral in estimate (1.6) is the volume growth of  $({\R}^n, g)$ when $q
= 2n/ (n - 2)$. In order to obtain (1.5) and (1.6) for large $p$ and $q$, additional
conditions on
$K$ or $u$ are required. By using a novel version of the moving plane method, Chen-Lin (\cite{Chen-Lin.1}
\cite{Chen-Lin.2} and \cite{Chen-Lin.3}) and Lin \cite{Lin.2} examine, among other things,
slow decay of $u$ under the condition 
$$
0 \,< \,{{C_7}\over {|x|^{1 + \alpha} }}\, \le \,|\btd K (x)| \,\le \,{{C_8}\over {|x|^{1 +
\alpha} }} \mfor {\mbox{large}}  \ \ |x|
\leqno (1.8)
$$
and for some positive constants $\alpha$, $C_7$ and $C_8$.\bk
To gain better understanding on $u$, consider the case when $K$ is equal to a positive
constant, say $K = n (n - 2)/4$, outside a compact subset of ${\R}^n$. We express $u$ as an associated function on the cylinder $\R \times S^{n - 1}$ by letting
$$
v (s, \theta) = |x|^{{n - 2}\over 2} u (x)\,, \ \ \ \ {\mbox{where}} \ \ |x| = e^s \ \ 
{\mbox{and}}
\
\
\theta = x / |x| \in S^{n - 1}. \leqno (1.9)
$$
Then $v$ satisfies the equation
$$
{{\partial^2 v}\over {\partial s^2}} + \Delta_\theta\, v -  {{ (n - 2)^2}\over 4}  v 
+ K  \,v^{{n + 2}\over {n - 2}} = 0 \ \ \ \
{\mbox{in}} \ \ \R \times S^{n - 1}, \leqno (1.10)
$$
where $\Delta_\theta$ is the standard Laplacian on $S^{n - 1}$. Here $K$ is interpreted as a
function on $\R \times S^{n - 1}$ such that $(s, \theta) \mapsto K (e^s, \theta)$ for $s \in \R$ and $\theta \in S^{n - 1}$.  By a
result of Caffarelli, Gidas and Spruck   
\cite{Caffarelli-Gidas-Spruck}, with improvements by Korevaar, Mazzeo, 
Pacard and  Schoen \cite {K-M-P-S},  either $g$  can be            
realized as a smooth metric on $S^n$ (in this case $u$ is said to have {\it fast decay}), or 
$$
v (s, \theta) = v_\varepsilon (s + T) \,[1 + O (e^{-\kappa s})] \mfor \ \ {\mbox{large}} \ \ s\,, 
\ \ \theta \in S^{n - 1} \leqno (1.11) 
$$
and for some  constants $\kappa > 0$ and $T \in \R$. Here $v_\varepsilon$, $\varepsilon \in (0, \, \, [(n - 2)/n]^{(n - 2)/4}]$,  is one of a one-parameter family of positive 
solutions of the O.D.E.
$$
v'' - {{(n - 2)^2}\over 4} \,v + {{n (n - 2)}\over 4}\,
v^{{n + 2}\over {n - 2}} = 0 \ \ \ \ {\mbox{in}} \ \ {\R}\,, \leqno (1.12)
$$
and $\varepsilon = {\mathop{\min}\limits_{t \in \R}} \ v (t)$ is referred as the necksize of the solution \cite {K-M-P-S}.   
As O.D.E. (1.12) is autonomous, 
$|v'_\varepsilon|$ is bounded in $\R$. Furthermore, the Pohozaev number 
$$
P (u) = \lim_{r \to + \infty} P (u, r) \ \ \ {\mbox{where}} \ \ \ P (u, r) = {{n - 2}\over {2n}} 
\int_{B_o (r)} x
\cdot
\btd K (x) \,u^{{2n}\over {n - 2}} (x)\, dx \leqno (1.13)
$$
is a negative number \cite{K-M-P-S}. When $K$ may not be a constant outside a compact subset
of ${\R}^n$, we have the following results.\\[0.2in] 
{\bf Theorem A.} \ \ {\it Let $u$ be a positive smooth solution of equation (1.1) with
condition (1.2), and $v$ given by (1.9). Assume that there exist positive constants
$C_9$ and $C_{10}$ such that}
$$
\int_{S^{n - 1}} \left( {{\partial v}\over {\partial s}} \right)^2 (s, \theta) \,d \theta 
\,\le\,
C_9 + C_{10}  
\int_{S^{n - 1}} v^2 (s, \theta)  \,d \theta \leqno (1.14)
$$
{\it for large $s$. If $P (u, r) \ge - \delta^2$ for large $r$ and for a positive constant
$\delta$, then }
$$
\int_{B_o (r)} u^{{2n}\over {n - 2}} \, dx \le C' \, \ln r \ \ \ \ and \ \ \ \ 
\int_{B_o (r)} |\btd u|^2 \, dx \le C'' \, \ln r \leqno (1.15)
$$
{\it for large $r$ and for some positive constants $C'$ and $C''$.}\\[0.2in]
{\bf Theorem B.} \ \ {\it Assumed that there exist positive constants $C_{11}$ and $C_{12}$ such that}
$$
\bigg\vert {{\partial v}\over {\partial s}} \bigg\vert  (s, \theta)  \le C_{11} + C_{12} \,v
(s, \theta) \ \ \ \ for \ \ large \ \ s \ \ and \ \ \theta \in S^{n - 1}. \leqno (1.16)
$$
{\it If $P (u, r) \ge - \delta^2$ for large $r$ and for a positive constant
$\delta$, then}
$$
\int_{S^{n - 1}} u^{{2n}\over {n - 2}} (r, \theta) \, d\theta \le C \,r^{-n} \leqno (1.17)
$$
{\it for large $s$ and for some positive constant $C$. Moreover, $u$ has slow decay.}\\[0.2in]
\hspace*{0.5in}We prove theorems A and B in section 4. Lower bounds on  
$P (u, r)$ are obtained in section 3, and examples are constructed in section 2. We use $c$, $C$, $C_1$, $C_2\,,...$ to denote positive constants, which may be different from section to section.

\vspace{0.5in}

{\bf \Large {\bf 2. \ \ Examples}}

\vspace{0.3in}

We begin with a construction of positive $C^2$-solution $u$ of equation (1.1) 
with $K$ bounded between two positive constants in ${\R}^n$, such that $u$ is unbounded from
above in ${\R}^n$ (and hence does not have slow decay), and the conformal metric $g$ is
complete. Throughout this note $n \ge 3$ is an integer. Let 
$$
\bar u \, (r, \lambda) = \alpha_n \, \left( {\lambda \over {\lambda^2 + r^2 }}
\right)^{(n - 2)/2} \mfor r \ge 0 \ \ {\mbox{and}} \ \ \lambda > 0\,, \leqno (2.1)
$$
where $\alpha_n = [n (n - 2)]^{(n - 2)/4},$ and 
$$
u_o (x) = \bar u \,( |x|\,, 1) = {{\alpha_n} \over { (1 + |x|^2)^{(n - 2)/2} }}
\mfor x \in {\R}^n. \leqno (2.2)
$$
Let $\{ \varepsilon_k \}_{k =
1}^\infty \subset (0, 1)$ be a sequence of decreasing numbers such that 
$$
\sum_{k = 1}^\infty \varepsilon_k = 1\,, \leqno (2.3)
$$
$\{ r_k\}_{k = 1}^\infty$  a sequence of positive numbers such that 
$r_1 \ge 1\,,$ $r_{k + 1} - r_k \ge 1$ for $k = 1, \,2,...,$  
and $\{ M_k \}$ a sequence of positive numbers such that $M_k \to + \infty$ as $k \to +\infty\,.$ 
For  
$
x^{1, k} := (r_k, \,0,..., 0) \in {\R}^n,$ $k = 1, \, 2,...,$ 
there exist positive numbers $\lambda_k\,,$ $k = 1, \, 2,...,$ such that 
$$
u_k (x) := \bar u \,(|x - x^{1, k}|, \,\lambda_k) \mfor x \in {\R}^n \leqno (2.4)
$$
satisfies 
$$
\Delta u_k + u_k^{{n + 2}\over {n - 2}} = 0\ \ \ \ {\mbox{in}} \ \ {\R}^n, \leqno (2.5)
$$
$$
u_k (x) \le \varepsilon_k \,u_o (x)  \ \ \ {\mbox{and}} \ \ \ |\btd u_k (x)\,| <
\varepsilon_k \mfor \ |x - x^{1, k}| \ge {1\over 4}\,,  \ \ {\mbox{and}} \leqno (2.6)
$$
$$
u_k (x^{1, k}) = \alpha_n \, \lambda_k^{(2 - n)/2}   \ge M_k 
\leqno (2.7)
$$
for $k = 1, 2,....$ Using (2.3) and (2.6), it follows as in \cite{Taliaferro.1} that ${\mathop{\sum}\limits_{k = 0}^\infty} \,u_k$ converges uniformly on compact subsets of ${\R}^n$ to a positive
$C^2$-function. For a positive number
$b$, let 
$$
\tilde u_b (x) =  (|x|^2 + b^2)^{(2 - n)/4}  \mfor x \in {\R}^n. \leqno (2.8)
$$
We have
$$
\Delta \,\tilde u_b + K_b \,\tilde u_b^{{n + 2}\over {n - 2}} = 0 \ \ \ \ {\mbox{in}} \ \ {\R}^n, 
\leqno (2.9)
$$
where 
$$
K_b (x) = {{n (n - 2)}\over 2} \left(1 - {{n + 2}\over {2n}} {{|x|^2}\over {|x|^2 + b^2}}
\right) 
\mfor x \in {\R}^n. \leqno (2.10)
$$
In particular
$$
{{n (n - 2)^2}\over {4n}} \le K_b (x) \le {{n (n - 2)}\over 2} \mfor x \in {\R}^n. \leqno
(2.11)
$$
Let 
$$
u (x) = \tilde u_b (x) + \sum_{k = 0}^\infty u_k (x) \mfor x \in {\R}^n.\leqno (2.12)
$$
It follows from (2.5), (2.9) and (2.11) that 
$$
- \Delta \,u (x) = \left[ K_b (x) \tilde u_b^{{n + 2}\over {n - 2}} (x) + \sum_{k = 0}^\infty 
u_k^{{n + 2}\over {n - 2}} (x) \right] \le {{n (n - 2)}\over 2} u^{{n + 2}\over {n - 2}} (x) 
  \leqno (2.13)
$$
for $x \in {\R}^n$. Assume that $x \in B_{x_{k'}} (1/4)$ for some positive integer $k'$. Using
(2.3), (2.6) and the inequality  
$
(a + b)^p \le 2^{p - 1} (a^p + b^p)  
$ 
for $a, \ b  \ge 0$ and $p \ge 1$,  we have
\begin{eqnarray*}
(2.14) \, \, u^{{n + 2}\over {n - 2}} (x) & = & \left[ \tilde u_b (x) + u_{k'} (x)  + \sum_{k \not= k'} u_k
(x) \right]^{{n + 2}\over {n - 2}}
 \le \left[ \tilde u_b (x) + u_{k'} (x) + 2 u_o (x) \right]^{{n + 2}\over {n - 2}}\\
& \ & \ \ \le \,c_1 \, \left[ \tilde u_b^{{n + 2}\over {n - 2}} (x) + u_{k'}^{{n + 2}\over {n - 2}}
(x)   + u_o^{{n + 2}\over {n - 2}} (x)  \right] 
\,\le \, - c_2 \,\Delta u (x)\,,
\end{eqnarray*}
where $c_1$ and $c_2$ are positive constants depending on $n$ only. Similar estimate holds for
$x
\not\in B_{x_{k'}} (1/4)$ for $k' = 1, \,2,...$, if we choose $c_2$ to be large enough, which 
depends on $n$ only. Thus $u$ satisfies the equation
$\Delta u + K u^{(n + 2)/(n - 2)} = 0$ in ${\R}^n,$ 
where 
$$
K (x) = [- \Delta u (x)] \,[u^{{n + 2}\over {n - 2}} (x)]^{-1} \mfor x \in {\R}^n  
$$
is a continuous function which is bounded in ${\R}^n$ between two positive constants by (2.13)
and (2.14). (2.7) shows that 
$u$ is  not bounded from above in ${\R}^n$. The conformal metric $u^{4/(n - 2)} g_o$ is
complete because
$$
u^{4/(n - 2)} (x) \ge \tilde u_b^{4/(n - 2)} (x) \ge (1/2) |x|^{-2} \leqno (2.15)
$$
for large $|x|$. 
Let 
$$
V_n := \omega_n \int_0^\infty \left( {{ \lambda \sqrt{n (n - 2)} }\over {\lambda^2 + r^2 }}
\right)^n r^{n - 1} \, dr = \omega_n \int_0^\infty \left( {{  \sqrt{n (n - 2)} }\over {1 + t^2
}}
\right)^n t^{n - 1} \, dt  \leqno (2.16)
$$
for $\lambda > 0$, where $t = \lambda^{-1}\, r$ and $\omega_n$ is the volume of the unit sphere in ${\R}^n$. 
By choosing $r_k$ suitably far from each other, together with (2.16) and the fact that the first integral in
(2.16) concentrates more on a neighborhood of $0$ for smaller $\lambda$, we have
$$
\int_{B_o (r)} u^{{2n}\over {n - 2}} (x) \, dx \le C_2 \, \ln r \leqno (2.17)
$$
for large $r$ and for a positive constant $C_2$. \bk
Next, given a function $\phi : [0, \infty) \to [0, \infty)$, we construct a positive $C^2$-solution $u$ of equation (1.1) with $K$ bounded between two positive constants in ${\R}^n,$  such that the conformal metric $g = u^{4/(n - 2)} g_o$ is complete and 
$$
\int_{B_o (r)} u^{{2n}\over {n - 2}} (x) \, dx \,\ge \,\phi (r) \mfor r > 2\,. \leqno (2.18)
$$
Without loss of generality, we may
assume that $\phi$ is increasing and $\phi (0) \ge 10 \, V_n$. For
$k = 1, \, 2,...$, let
$N_k$ be a positive integer such that 
$$
N_k \ge  2 \,V_n^{-1}  \,\phi \,(k + 2) \mfor k = 1, \, 2,.... \leqno (2.19)
$$
Let $\{ \epsilon_k \}_{k =
1}^\infty \subset (0, 1)$ be a sequence of decreasing numbers such that 
$$
\sum_{k = 1}^\infty N_k\,\epsilon_k \,\le \,1\,. \leqno (2.20)
$$
Let $\theta_k = 2
\pi/ N_k$. Let
$$
x_{k, j} = \left(k \, \sin \, (j \,\theta_k)\,, \  
k \cos \,(j \,\theta_k)\,, \ 0\,,..., \ 0 \right) \ \in {\R}^n \mfor j = 1, \, 2,..., N_k\,,
\leqno (2.21)
$$
and  
$$
u_{k, j} (x) =  \bar u \,(|x - x_{k, j} \,|, \,\lambda_k) \mfor x \in {\R}^n \ \ \ \
{\mbox{and}} \ \ j = 1, \ 2,..., N_k\,.
\leqno (2.22)
$$
We choose $\lambda_k$ to be small so that 
$$
u_{k, j} (x) \le \epsilon_k \,u_o (x)  \ \ \ {\mbox{and}} \ \ \ |\btd u_{k, j} (x)\,| <
\epsilon_k \mfor \ |x - x_{k, j}\,|                    \ge  \pi / (10 N_k) \,, \leqno (2.23)
$$
and 
$$ 
\int_{B_{x_{k, j}} (\pi/ (10 N_k))}  u_{k, j}^{{2n}\over {n - 2}} (x) \, dx \ge {{V_n}\over
2} \mfor j = 1, \, 2,..., N_k\,, \leqno (2.24)
$$
where $B_{x_{k, j}} (\pi/ (10 N_k))$ is the ball with center at $x_{k, j}$ and radius equal to $\pi/ (10 N_k)$. 
(2.24) is possible because, when $\lambda$ is smaller, the first integral in (2.16)
concentrates more on a neighborhood of the origin. It follows from (2.20) and (2.23) that the series ${\mathop{\sum}\limits_{k = 1}^\infty} \, {\mathop{\sum}\limits_{j = 1}^{N_k}} \, u_{k, j}$ converges 
uniformly on compact subsets of ${\R}^n$ to a positive $C^2$-function.  
Let 
$$
u = {\tilde u}_b + u_o + {\mathop{\sum}\limits_{k = 1}^\infty} \, {\mathop{\sum}\limits_{j = 1}^{N_k}} \, u_{k, j} \ \ \ \ {\mbox{in}} \ \ {\R}^n.
$$ 
As above, we have
$\Delta u + K u^{(n + 2)/(n - 2)} = 0$ in ${\R}^n,$ 
where $K$ is a continuous function on ${\R}^n$ that is bounded between two positive
constants. For any $r > 2$, let $k$ be the integer such that $k + 1 \le r < k +
2$. By (2.19) we have
\begin{eqnarray*}
\phi (r) & \le & \phi (k + 2) \,\le\,  {{V_n N_k}\over 2} \, \le \, \sum_{j =1}^{N_k} \int_{B_{x_{k, j}}
(\pi/ (10 N_k))}  u_{k, j}^{{2n}\over {n - 2}} (x) \, dx\\
& \le & \int_{B_o (k + 1)}  u^{{2n}\over {n - 2}} (x) \, dx \,\le \,\int_{B_o (r)} 
u^{{2n}\over {n - 2}} (x) \, dx\,. 
\end{eqnarray*}


\vspace{0.4in}

{\bf \Large {\bf 3. \ \ Estimates on $P (u, r)$}}

\vspace{0.3in}

Let $P (u, r)$ be given by (1.13) in the introduction. The Pohozaev identity (see, for example, \cite{Ding-Ni}) states that  
$$
P (u, r) =  \int_{S_r} \left[ r \left( {{\partial u}\over {\partial r}} \right)^2
 - {r\over 2} |\btd u|^2 +  {{n - 2}\over {2n}}\, r  K u^{{2n}\over {n -
2}}  + {{n - 2}\over 2} \,u \,{{\partial
u}\over {\partial r}} \right] \,dS  \leqno (3.1)
$$
for $r > 0$, where $S_r = \partial B_o (r)$ is the sphere of radius
$r$. \\[0.2in]
{\bf Theorem 3.2.} \ \ {\it Let $u$ be a positive smooth solution of equation (1.1)
with condition (1.2). Assume that
$u$ is bounded from above in
${\R}^n$ and}
$$
{{\partial K}\over {\partial r}} (x) \, \ge \, -\,{{C_1}\over { |x|^{(n + 2)/2} \,(\ln \,
|x|)^{1 + \epsilon} }} \leqno (3.3)
$$
{\it for large $|x|$ and for some positive constants $C_1$ and $\epsilon$. 
Then $P (u, r) \ge - \delta^2$ for large $r$ and for a positive constant $\delta$.}\\[0.1in]
{\bf Proof.} \ \ Fixing a large number $R$ and using (1.4) we have
\begin{eqnarray*}
& \ & \int_{ B_o ( (m + 1) R ) \setminus B_o (m R)} r {{\partial K}\over {\partial r}} (x)
u^{{2n}\over {n - 2}} (x) \, dx\\ 
& \ge & - {{C_1}\over { (m R)^{{n}\over 2} [\ln \, (m
R)]^{1 + \epsilon} }} \int_{ B_o ( (m + 1) R ) \setminus B_o (m R)} u^{{2n}\over {n - 2}} (x)
\, dx\\
& \ge & - {{C_2}\over { (m R)^{{n}\over 2} [\ln \, (m
R)]^{1 + \epsilon} }} \int_{ B_o ( (m + 1) R ) \setminus B_o (m R)} u^{{n + 2}\over {n - 2}}
(x)
\, dx\\
& \ge & - \,{{C_3 [ (m + 1) R]^{{n - 2}\over 2} }\over { (m R)^{{n}\over 2} [\ln \, (m R)]^{1
+
\epsilon} }} \, \ge \, - \,{{C_4}\over {m \,(\ln \,m)^{1 + \epsilon} }}
\end{eqnarray*}
for any positive integer $m$ larger than $1\,,$ where $r = |x|$. Here $C_2$, $C_3$ and $C_4$
are positive constants. As the series
$$
\sum_{m = 2}^{\infty} {1 \over {m \,(\ln m)^{1 + \epsilon} }}
$$
converges, we conclude that there exists a positive constant $\delta$ such that $P (u, r) \ge - \,\delta^2$ for large $r$. \qedwh
{\bf Theorem 3.4.} \ \ {\it Let $u$ be a positive smooth solution of equation (1.1)
with condition (1.2). 
Assume that there exists a positive constant $c$ such that}
$$
{{\partial K}\over {\partial r}} (r, \theta) \,\ge \, - {c\over {r^2}} \ \ \ \ for \ \ large \ \ r \ \ and \ \ \theta \in S^{n - 1}. \leqno (3.5)
$$
{\it  If there exist positive constants $C$ and
$\lambda \in (0, 1)$ such that}
$$
\int_{S^{n - 1}} \left( {{\partial v}\over {\partial s}} \right)^{{2n}\over {n - 2}} 
(s, \theta) \, d \theta \,\le C \,\, e^{\lambda \,s} \leqno (3.6)
$$
{\it for large $s$, then  $P (u, r) \ge - \,\delta^2$ for large $r$ and for a positive constant $\delta$.}

\pagebreak

{\bf Proof.} \ \ For a positive number $\varepsilon > 0$ such that 
$\varepsilon + \lambda < 1$, using Young's inequality we have
\begin{eqnarray*}
& \ & {d\over {dr}} \left( \int_{S_r} r^\varepsilon u^{{2n}\over {n - 2}} (r, \theta)\, dS \right)
\,= \,{d\over {dr}} \left( \int_{S^{n - 1}} r^{n - 1 + \varepsilon} u^{{2n}\over {n - 2}} (r, \theta)\, d\theta \right)\\
& = & {{ n - 1 + \varepsilon}\over r} \int_{S^{n - 1}} r^{n - 1 + \varepsilon} 
u^{{2n}\over {n - 2}} (r, \theta)\, d\theta + {{2n}\over {n - 2}} \int_{S^{n - 1}} u^{{n +
2}\over {n - 2}} (r, \theta) {{\partial u}\over {\partial r}}  (r, \theta) \, r^{n - 1 +
\epsilon} \, d\theta\\ 
& = & {{ - 1 + \varepsilon}\over r} \int_{S^{n - 1}} r^{n - 1 +
\varepsilon} u^{{2n}\over {n - 2}} (r, \theta)\, d\theta\\
& \ & \ \ \ \ \ \, + \,{{2n}\over {n - 2}} \int_{S^{n -
1}} u^{{n + 2}\over {n - 2}}  (r, \theta)
\left[ {{\partial u}\over {\partial r}} + {{n - 2}\over 2} {u \over r} \right]  (r, \theta) \,
r^{n - 1 + \epsilon} \, d\theta\\ 
& \le & {{C_5}\over {r^{2 - \varepsilon}}} \int_{S^{n - 1}} \left\{ r^{n\over 2}
\left[ {{\partial u}\over {\partial r}} + {{n - 2}\over 2} {u \over r} \right]  (r, \theta)
\right\}^{{2n}\over {n - 2}} d\theta \\
& = & \, {{C_5}\over {r^{2 - \varepsilon}}} \int_{S^{n - 1}} 
\left( {{\partial v}\over {\partial s}} \right)^{{2n}\over {n - 2}} 
(s, \theta) \, d \theta
\, \,\le \, \,{{C_6}\over {r^{ 2 - \lambda - \varepsilon}}}
\end{eqnarray*}
for large $r$, where $r = e^s$ and $C_5$ and $C_6$ are positive constants. It follows that there exists a
positive constant $C_7$ such that 
$$
\int_{S_r} r^\varepsilon u^{{2n}\over {n - 2}} \, dS  \le C_7
\ \ \ \ {\mbox{or}} \ \ \ \  
\int_{S_r} u^{{2n}\over {n - 2}} \, dS  \le C_7 \,r^{- \varepsilon} \leqno (3.7)
$$
for large $r$. For a fixed large number $R_o\,,$ we have
$$
{{2n}\over {n - 2}} \,P (u, R) = \int_{B_o (R)} r {{\partial K}\over {\partial r}} u^{{2n}\over {n - 2}} dx \ge - C_8  - C_9
\int_{R_o}^R r^{-1} \int_{S_r}  u^{{2n}\over {n - 2}} \, dS \, dr \ge - C_{10}
$$
for large $R$ with $R_o < R$. Here $C_8$, $C_9$ and $C_{10}$ 
are positive constants.\qedwh


\vspace{0.3in}

{\bf \Large {\bf 4. \ \ Proofs of Theorem A and B}}

\vspace{0.3in}

{\bf Proof of Theorem A.} \ \ Let 
$$
w (s) = {1\over 2} \int_{S^{n - 1}} v^2 (s, \theta)  \,d \theta \ \ \ \ {\mbox{for}} \ \ s \in
\R\,, \leqno (4.1)
$$
where $v$ is defined in (1.9). Using equation (1.10) we have 
\begin{eqnarray*}
(4.2) \ \ \ w'' (s) & = & \int_{S^{n - 1}} \left( {{\partial v}\over {\partial s}}
\right)^2 (s,
\theta)
\,d
\theta + \int_{S^{n - 1}} |\btd_\theta v  \,(s, \theta) |^2 \, d \theta
\ \ \ \ \ \ \\ & \ & \ \ \ + \left( {{n - 2}\over
2} \right)^2 \int_{S^{n - 1}} v^2 (s, \theta)  \,d \theta - \int_{S^{n - 1}} K (e^s, \theta) \,v^{{2n}\over {n
- 2}}  (s,
\theta) 
\,d \theta  \ \ \ \  \ \ \ \  \ \  \  
\end{eqnarray*}
for $s \in \R\,.$ The Pohozaev identity can be expressed as
\begin{eqnarray*}
(4.3)  \ \ 2 P (u, e^s) & = & \int_{S^{n - 1}} \left( {{\partial v}\over {\partial
s}}
\right)^2 (s,
\theta)
\,d
\theta - \int_{S^{n - 1}} |\btd_\theta v  \,(s, \theta) |^2 \, d \theta 
\ \ \ \ \ \ \\
& \ & \ \ 
- \left( {{n - 2}\over
2} \right)^2 \int_{S^{n - 1}} v^2 (s, \theta)  \,d \theta + {{n - 2}\over n} \int_{S^{n - 1}}
K (e^s, \theta) \,v^{{2n}\over {n - 2}}  (s,
\theta) 
\,d \theta  \ \ \ \  \ \ \ \  \ \ \ \  
\end{eqnarray*}
for $s \in \R$ \cite{Cheung-Leung.1}. It follows from (1.2), (1.4), (4.2) and (4.3) that
$$
w'' (s) \le C_1 + C_2 \int_{S^{n - 1}} v^2
(s,
\theta)
\,d
\theta - {{2 \, a^2} \over n} \int_{S^{n - 1}} v^{{2n}\over {n
- 2}}  (s, \theta)  \,d \theta \leqno (4.4)
$$
for large $s$, where $C_1$ and $C_2$ are positive constants. Applying Young's inequality we
obtain
$$
 - \,{{2\, a^2}\over n} \int_{S^{n - 1}}  v^{{2n}\over {n
- 2}}  (s, \theta)  \,d \theta \le C_3 - C_4 \int_{S^{n - 1}} v^2
(s,
\theta)
\,d
\theta \leqno (4.5)
$$
for large $s$, where $C_3$ and $C_4$ are positive constants. Furthermore, by choosing $C_3$
to be large, we can take $C_4$ to be large as well. Hence there exists a positive constant
$C_5$ such that 
$$
w'' (s) \le C_5 - w (s) \ \ \ \ {\mbox{for \ \ large}} \ \ s\,. \leqno (4.6)
$$
From (4.6) it is easy to see that $w (s)$ is uniformly bounded from above for large $s$. To prove this assertion, 
assume that there is a large $s'$ such that $w (s') \ge C_5 + 1$. (4.6) implies that $w'' (s')
\le -1$. Let $s_o$ be a number larger than $s'$ such that 
$w (s_o) < C_5 + 1$ and $w' (s_o) \le 0\,.$ 
If $w (s) < C_5 + 1$ for all $s > s_o$, then we are done. Assume that $s_1$ is the smallest
number larger than $s_o$ such that $w (s_1) = C_5 + 1$. We claim that 
$$
D := w' (s_1) < 2 \,(C_5 +
1)\,. \leqno (4.7)
$$
Let $\bar s \in (s_o, s_1)$ be the largest number such that $w' (\bar s) = D/2$. As
$w''
\le C_5$  on $(s_o\,,
s_1)$, we have  
$s_1 - \bar s \ge D/(2C_5)\,.$ 
On the other hand, $w' \ge D/2$ on $(\bar s, s_1)$. Therefore we have
$$
C_5 + 1 \,\ge \,w (s_1) - w (\bar s) \,\ge\, {D\over {2C_5}} \cdot {D\over 2} \,\,\Rightarrow
\,\, D^2
\,\le \,4 C_5 (C_5 + 1)\,. 
$$
Hence we have (4.7). From $s_1$, $w (s)$ can become no larger than $(C_5 + 1) + [2 (C_5 +
1)]^2$ before
$w' (s)$ becomes negative again. Hence we conclude that $w (s)$ is uniformly bounded from
above for large
$s$.\bk 
From Pohozaev identity (3.1) we obtain
$$
\int_{S_r} r |\btd u|^2 \, dS = 
2  \int_{S_r} \left[ r \left( {{\partial u}\over {\partial r}} \right)^2
+  {{n - 2}\over {2n}} r  K u^{{2n}\over {n -
2}}  + {{n - 2}\over 2} \,u \,{{\partial
u}\over {\partial r}} \right] \,dS - 2 P (u, r) \leqno (4.8)
$$
for $r > 0$. We have
$$
\int_{S_r}  r \left( {{\partial u}\over {\partial r}} \right)^2 \, dS 
= \int_{S_r}  r \left[ {{\partial u}\over {\partial r}} + {{n - 2}\over 2} {u \over
r} \right]^2 dS - (n - 2) \int_{S_r} u {{\partial u}\over {\partial r}}  \, dS - 
\left( {{n - 2}\over 2}  \right)^2 \int_{S_r} {{u^2}\over r} \, dS \leqno (4.9)
$$
for $r > 0\,.$ Using (1.14) and the fact that $w$ is bounded from above we obtain
\begin{eqnarray*}
(4.10) \ \ \ \  - \int_{S_r} u {{\partial u}\over {\partial r}}  \, dS & = & \, - \int_{S_r} u \left[  
{{\partial u}\over {\partial r}} + {{n - 2}\over 2} {u \over
r} \right] \,dS + {{n - 2}\over 2} \int_{S_r} {{u^2}\over r} \, dS\\
& \le & \int_{S_r} r\left[  
{{\partial u}\over {\partial r}} + {{n - 2}\over 2} {u \over
r} \right]^2 \,dS + {n\over 2} \int_{S_r} {{u^2}\over r} \, dS\\
& = & \int_{S^{n - 1}} \left( {{\partial v}\over {\partial s}} \right)^2 (s, \theta) 
\,d \theta 
+ {n\over 2} \int_{S^{n - 1}} v^2 (s, \theta)  \,d \theta \le C_6 \ \ \ \ \ \ \ \ 
\end{eqnarray*}
for large $r$ and for a positive constant $C_6\,,$ where $r = e^s$. It follows from (4.8), (4.9) and (4.10) that 
$$
\int_{S_r} |\btd u|^2 \, dS \,\le \, {{C_7}\over r} 
+ {{n - 2}\over {n}} \int_{S_r} K u^{{2n}\over {n - 2}} \, dS  
$$
for large $r$, where $C_7$ is a positive constant. Therefore we obtain 
$$
\int_{B_o (r)} |\btd u|^2 \, dx \,\le \,C_8 \ln r 
+ {{n - 2}\over {n}} \int_{B_o (r)} K u^{{2n}\over {n - 2}} \, dx \leqno (4.11)
$$
for large $r$ and for a positive constant $C_8 \ge C_7$. On the other hand we have
\begin{eqnarray*}
\int_{B_o (r)} K u^{{2n}\over {n - 2}} \, dx & = & \int_{B_o (r)} u \,( - \Delta u ) \, dx 
= \int_{B_o (r)} |\btd u|^2 \, dx - \int_{S_r} u {{\partial u}\over {\partial r}} \, dS\\
& \le & C_8 \ln r + {{n - 2}\over {n}} \int_{B_o (r)} K u^{{2n}\over {n - 2}} \, dx + C_6
\end{eqnarray*}
for large $r$, where we use (4.10). Hence there exists a positive constant $C_9$ such that 
$$
\int_{B_o (r)} K u^{{2n}\over {n - 2}} \, dx \le C_9 \, \ln r
$$
for large $r$. If $u \in L^{2n/(n - 2)} ({\R}^n)$, then clearly we have the first inequality in (1.15). Assume that $u \not\in L^{2n/(n - 2)} ({\R}^n)$. Using (1.2) we have 
$$
\int_{B_o (r)} u^{{2n}\over {n - 2}} \, dx \, \le \, {2\over {a^2}} \, \int_{B_o (r)} K u^{{2n}\over {n - 2}} \, dx \le {{2 \,C_9}\over {a^2}} \, \ln r
$$
for large $r$. Hence we have the first inequality in
(1.15). The second inequality follows from (4.11).\qedwh
{\bf Proof of Theorem B.} \ \  From the proof of theorem A we have
$$
\int_{S_r} {{u^2 (x)}\over r} \, dS = \int_{S^{n - 1}} v^2 (s, \theta) \, d\theta = 2 \,w (s) \le
C_{10}
\leqno (4.12)
$$
for large $r$, where $r = |x| = e^s$ and $C_{10}$ is a positive constant. By using (1.16) and (4.12) we also have
$$
\int_{S_r} r \left[ {{\partial u}\over {\partial r}} + {{n
- 2}\over 2} {u\over r} \right]^2 dS = \int_{S^{n - 1}} r^n 
\left[ {{\partial u}\over {\partial r}} + {{n - 2}\over 2} {u\over r} \right]^2 d\theta =
\int_{S^{n - 1}} \left( {{\partial v}\over {\partial s}}\right)^2 d\theta  \le C_{11}
\leqno (4.13)
$$
for large $r$, where $C_{11}$ is a positive constant. 
It follows from Pohozaev identity (3.1) that
\begin{eqnarray*}
(4.14)\,\int_{S_r} r | \btd u |^2 \,dS & \le & C_{12} + {{n - 2}\over {n}} \int_{S_r} r K
u^{{2n}\over {n - 2}}
\, dS + 2 \int_{S_r} r \left[ {{\partial u}\over {\partial r}} + {{n
- 2}\over 2} {u\over r} \right]^2 dS\\
& \ & \ \ \ \ \  \, - \, (n - 2)
 \int_{S_r} u \left[ {{\partial u}\over {\partial r}} + {{n
- 2}\over 2} {u\over r} \right] \,dS\\
& \le & C_{13} + {{n - 2}\over {n}} \int_{S_r} r K u^{{2n}\over {n - 2}} \, dS 
+ C_{14} \int_{S_r} r \left[ {{\partial u}\over {\partial r}} + {{n - 2}\over 2} {u\over r}
\right]^2 dS \\
& \ & \ \ \ \ \ \, + \, C_{15}
\int_{S_r} {{u^2}\over r} \, dS \, \le \, C_{16} + {{n - 2}\over {n}} \int_{S_r} r K
u^{{2n}\over {n - 2}} \, dS \ \ 
\end{eqnarray*}
for large $r$, where we use (4.12) and (4.13). Here $C_{12}$, $C_{13},$ $C_{14},$ $C_{15}$
and $C_{16}$ are positive constants.  
(4.14) implies 
that there exists a positive constant
$C_{17}$ such that  
$$
\int_{B_o (R)} r | \btd u |^2 \,dx \,\le \,C_{17} \, R + {{n - 2}\over {n}} \int_{B_o (R)} r K
u^{{2n}\over {n - 2}} \, dx \leqno (4.15)
$$
for large $R$. We have 
\begin{eqnarray*}
(4.16) \ \ \ \ \ \ \ \  & \ & \int_{B_o (R)} r K u^{{2n}\over {n - 2}} \, dx 
= \int_{B_o (R)} (ru) (-
\Delta u) \, dx\\ 
& = & \int_{B_o (R)} r | \btd u |^2 \, dx + \int_{B_o (R)} u
{{\partial u}\over {\partial r}} \, dx - R \int_{S_R} u {{\partial u}\over {\partial r}}
\, dS  \ \ \ \ \ \ \  \ \ \ \ \ \ \ \ \ \ \ \ \ \ \ 
\end{eqnarray*}
for $R > 0$. Using (4.13) we obtain 
\begin{eqnarray*}
(4.17) \ \ \ \ \int_{B_o (R)} u\,
{{\partial u}\over {\partial r}} \, dx & = & \int_{B_o (R)} u
\left[ {{\partial u}\over {\partial r}} + {{n - 2}\over 2} {u\over r} \right] \, dx - {{n -
2}\over 2} \int_{B_o (R)} {{u^2}\over r} \, dx \ \ \ \ \ \ \  \ \ \\
& \le & C_{18} \int_{B_o (R)} r 
\left[ {{\partial u}\over {\partial r}} + {{n - 2}\over 2} {u\over r} \right]^2 \, dx
\, \le \, C_{19} \, R
\end{eqnarray*}
for large $R$, where $C_{18}$ and $C_{19}$ are positive constants.  
As in (4.10) we have 
$$
\bigg\vert R \int_{S_R} u {{\partial u}\over {\partial r}}\, dS \bigg\vert \, \le \, C_{20} \,R \leqno
(4.18)
$$
for large $R\,,$ where $C_{20}$ is a positive constant. From (4.15), (4.16), (4.17) and (4.18) we obtain 
$$
{2\over n} \int_{B_o (R)} r K u^{{2n}\over {n - 2}} \, dx  \le C_{21} \,R \leqno (4.19)
$$
for large $R$, where $C_{21}$ is a positive constant. Using (1.16) we have
\begin{eqnarray*}
& \ & {d\over {dr}} \left( \int_{S_r} r^2 u^{{2n}\over {n - 2}} \, dS \right) 
= {d\over {dr}} \left( \int_{S^{n - 1}} r^{n + 1} u^{{2n}\over {n - 2}} \, d\theta \right)\\ 
& = & (n + 1) \int_{S^{n - 1}} r^n u^{{2n}\over {n - 2}} \, d\theta 
+ {{2n}\over {n - 2}} \int_{S^{n - 1}} r^{n + 1} u^{{n + 2}\over {n - 2}} {{\partial
u}\over {\partial r}} \, d\theta\\
& = &  \int_{S^{n - 1}} r^n u^{{2n}\over {n - 2}} \, d\theta 
+ {{2n}\over {n - 2}} \int_{S^{n - 1}} r^{n + 1} u^{{n + 2}\over {n - 2}} \left[ {{\partial
u}\over {\partial r}} + {{n - 2}\over 2} {u\over r} \right] \, d\theta\\
& = & \int_{S^{n - 1}} r^n u^{{2n}\over {n - 2}} \, d\theta + 
{{2n}\over {n - 2}} \int_{S^{n - 1}} \left(r^{{n + 2}\over 2} u^{{n + 2}\over {n - 2}}
\right) \left\{ r^{n\over 2} \left[ {{\partial u}\over {\partial r}} + {{n - 2}\over 2}
{u\over r}
\right] \right\} \, d\theta\\
& \le & C_{22} \int_{S^{n - 1}} r^n u^{{2n}\over {n - 2}} \, d\theta 
+  C_{23} \int_{S^{n - 1}} \left\{ r^{n\over 2} \left[ {{\partial u}\over {\partial r}} + {{n
- 2}\over 2} {u\over r}
\right] \right\}^{{2n}\over {n - 2}} \, d\theta\\
& = &  C_{22} \int_{S^{n - 1}} r^n u^{{2n}\over {n - 2}} \, d\theta  +  C_{23} \int_{S^{n -
1}}
\bigg\vert {{\partial v}\over {\partial s}} \bigg\vert^{{2n}\over {n - 2}} d\theta
\, \le \,  C_{24} \int_{S^{n - 1}} r^n u^{{2n}\over {n - 2}} \, d\theta  +  C_{25}
\end{eqnarray*}
for large $r$, and for some positive constants $C_{22}$, $C_{23}$, $C_{24}$ and $C_{25}$, where 
$r = e^s$. 
Hence 
\begin{eqnarray*}
(4.20) \ \ \ \ \ \  \int_{S_R} R^2 u^{{2n}\over {n - 2}}
\, dS & = & \int_0^R  
 \left( \int_{S_t} t^2 u^{{2n}\over {n - 2}} \, dS \right)' dt
\,\le \, C_{26} + \int_{r_o}^R 
 \left( \int_{S_t} t^2 u^{{2n}\over {n - 2}} \, dS \right)' dt\\
& \le & C_{26} \,+\, C_{24}\int_{r_o}^R \int_{S_t} t u^{{2n}\over {n - 2}} \, dS \,dt + C_{25} (R
- r_o)\\ & \le & C_{27} \,R \,+\, C_{28} \int_{B_o (R)} r  u^{{2n}\over {n - 2}} \,dx \ \ \ \ \ 
\ \ \ \ \ \ \ \ \ \ \ \ \ \ \ \  
\end{eqnarray*}
for $R$ and $r_o$ large, with $R > r_o$. Here $C_{26}$, $C_{27}$ and $C_{28}$ are positive constants.   
Consider the case when $u \not\in L^{2n/(n - 2)} ({\R}^n)$. It follows from (1.2) and (4.19) that 
$$
\int_{B_o (R)} r  u^{{2n}\over {n - 2}} dx \,\le\, C_{29} \int_{B_o (R)} r  K u^{{2n}\over {n -
2}} dx 
\le C_{30} \,R \leqno (4.21)
$$
for large $R$ and for some positive constants $C_{29}$ and $C_{30}$. Clearly we have
$$
\int_{B_o (R)} r  u^{{2n}\over {n - 2}} dx \,\le\, C_{31} \,R \leqno (4.22)
$$
for large $R$ and for some positive constants $C_{31}$ if $u \in
L^{2n/(n - 2)} ({\R}^n)\,.$ From (4.20), (4.21) and (4.22) we have (1.17). By the results in \cite{Leung.5} (see also \cite{Cheung-Leung.1}), we obtain slow decay
(1.7)  as well.\qed.

\vspace*{5in}

{\large D}EPARTMENT OF {\large M}ATHEMATICS, 
{\large N}ATIONAL {\large U}NIVERSITY OF {\large S}INGAPORE, 
2 {\large S}CIENCE {\large D}RIVE 2, {\large S}INGAPORE 117543, 
{\large R}EPUBLIC OF {\large S}INGAPORE; \ \ \ \ 
{\tt matlmc@math.nus.edu.sg}

\end{document}